\documentclass[12pt,leqno]{article}
\usepackage[a4paper, margin=25.4mm]{geometry}
\usepackage{amssymb}
\usepackage{amsmath}
\usepackage{amsthm}
\usepackage{stmaryrd}
\usepackage{mathtools,mathrsfs}
\usepackage[pagebackref]{hyperref}
\hypersetup{citecolor=blue, linkcolor=blue, colorlinks=true} 
\usepackage{booktabs}
\usepackage{enumerate}
\usepackage{tikz}
\usepackage{url}
\usepackage[T2A,T1]{fontenc}
\usepackage{authblk}  

\newtheorem{theorem}{Theorem}[section]
\newtheorem{proposition}[theorem]{Proposition}
\newtheorem{lemma}[theorem]{Lemma}
\newtheorem{corollary}[theorem]{Corollary}
\theoremstyle{definition}

\newtheorem{remark}[theorem]{Remark}

\renewcommand{\le}{\leqslant}

\newcommand{\coloneq}{\vcentcolon=}      

\newcommand{\eps}{\varepsilon}
\newcommand{\cH}{{\mathcal{H}}}
\newcommand{\cK}{{\mathcal{K}}}
\newcommand{\cT}{{\mathcal{T}}}

\newcommand{\R}{{\mathbb R}}
\newcommand{\Z}{{\mathbb Z}}

\def\Kop{\operatornamewithlimits{%
  \mathchoice{\vcenter{\hbox{\LARGE $\mathcal{K}$}}}
             {\vcenter{\hbox{\large $\mathcal{K}$}}}
             {\mathrm{\mathcal{K}}}
             {\mathrm{\mathcal{K}}}}}

\author[1]{S.P. Glasby}
\affil[1]{\small Center for the Mathematics of Symmetry and Computation,
  University of\newline Western Australia, Perth 6009, Australia\quad
  \href{mailto:Stephen.Glasby@uwa.edu.au}{Stephen.Glasby@uwa.edu.au}}

\author[2]{G.R. Paseman}
\affil[2]{\small Sheperd Systems, UC Berkeley\quad
  \href{mailto:sheperdsystems@gmail.com}{sheperdsystems@gmail.com}}

\title{Maximizing weighted sums of binomial coefficients using generalized continued fractions}

\date{\today}

\begin{document}
\maketitle

\begin{abstract}
  Let $m,r\in\Z$ and $\omega\in\R$ satisfy $0\leqslant r\leqslant m$ and
  $\omega\geqslant1$. Our main result is a generalized continued fraction
  for an expression  involving the partial binomial sum
  $s_m(r) = \sum_{i=0}^r\binom{m}{i}$.  We apply this to create new upper and
  lower bounds for $s_m(r)$ and thus for  $g_{\omega,m}(r)=\omega^{-r}s_m(r)$.
  We also bound an integer $r_0 \in \{0,1,\dots,m\}$
  such that $g_{\omega,m}(0)<\cdots<g_{\omega,m}(r_0-1)\leqslant
  g_{\omega,m}(r_0)$ and $g_{\omega,m}(r_0)>\cdots>g_{\omega,m}(m)$.
  For real $\omega\geqslant\sqrt3$ we prove that
  $r_0\in\{\lfloor\frac{m+2}{\omega+1}\rfloor,\lfloor\frac{m+2}{\omega+1}\rfloor+1\}$,
  and also $r_0 =\lfloor\frac{m+2}{\omega+1}\rfloor$ for
  $\omega\in\{3,4,\dots\}$ or $\omega=2$ and $3\nmid m$.
  \vskip2mm\noindent 
  {\bf Keywords:} partial sum, binomial coefficients, continued fraction,
      bounds
  \vskip2mm\noindent
  {\bf 2020 Mathematics Subject Classification:} 05A10, 11B65, 11Y65
\end{abstract}


\section{Introduction}

Given a real number $\omega\geqslant1$ and integers $m,r$ 
satisfying $0\leqslant r\leqslant m$, set
\begin{equation}\label{E1}
  s_m(r)\coloneq\sum_{i=0}^r\binom{m}{i}\qquad\textup{and}\qquad
  g(r)= g_{\omega,m}(r)\coloneq\omega^{-r}s_m(r),
\end{equation}
where
the binomial coefficient $\binom{m}{i}$ equals $\prod_{k=1}^i\frac{m-k+1}{k}$
for $i>0$ and $\binom{m}{0}=1$.
The weighted binomial sum $g_{\omega,m}(r)$ and the partial binomial 
sum $s_m(r)=g_{1,m}(r)$ appear in many formulas and inequalities,
e.g. the cumulative distribution function $2^{-m}s_m(r)$
of a binomial random variable with $p=q=\frac{1}{2}$ as in Remark~\ref{BE},
and the Gilbert-Varshamov bound~\cite[Theorem~5.2.6]{LX} for a code
$C\subseteq\{0,1\}^n$.
Partial sums of binomial coefficients are found in probability theory,
coding theory,
group theory, and elsewhere. As $s_m(r)$ cannot be computed exactly for most
values of $r$, it is desirable for certain applications
to find simple sharp upper and lower bounds for $s_m(r)$.
Our interest in bounding $2^{-r}s_m(r)$ was piqued in~\cite{GP} by an
application to Reed-Muller codes
$\textup{RM}(m,r)$, which are linear codes of dimension $s_m(r)$.

Our main result is a generalized continued fraction 
$a_0 + \Kop_{i=1}^{r} \frac{b_i}{a_i}$ (using Gauss' Kettenbruch notation)
for $Q\coloneq\frac{(r+1)}{s_m(r)}\binom{m}{r+1}$.
From this we derive useful approximations to $Q, 2+\frac{Q}{r+1}$,
and $s_m(r)$, and with these find a maximizing input $r_0$ for $g_{\omega,m}(r)$.

The \emph{$j^{\rm th}$ tail} of the generalized continued fraction
$\Kop_{i=1}^{r} \frac{b_i}{a_i}$ is denoted
by $\cT_j$ where
\begin{equation}\label{E:K}
  \cT_j\coloneq\Kop_{i=j}^{r} \frac{b_i}{a_i}=
  \cfrac{b_j}{a_j +\cfrac{b_{j+1}}{a_{j+1} +\cfrac{b_{j+2}}{\llap{\ensuremath{\ddots}} \raisebox{-0.8em}{\ensuremath{a_{r-1} + \cfrac{b_r}{a_r}}}}}}
  =\frac{b_j}{a_j+\cT_{j+1}}  \qquad\qquad\textup{and $1\leqslant j\leqslant r$.}
\end{equation}
If $\cT_j=\frac{B_j}{A_j}$, then $\cT_j=\frac{b_j}{a_j+\cT_{j+1}}$ shows
$b_jA_j-a_jB_j=\cT_{j+1}B_j$. By convention $\cT_{r+1}=0$.

It follows from $\binom{m}{r-i}=\binom{m}{r}\prod_{k=1}^i\frac{r-k+1}{m-r+k}$
that $x^i\binom{m}{r}\leqslant\binom{m}{r-i}\le y^i\binom{m}{r}$ for $0\leqslant i\leqslant r$
where $x\coloneq\frac{1}{m}$ and $y\coloneq\frac{r}{m-r+1}$. Hence
$\frac{1-x^{r+1}}{1-x}\binom{m}{r}\leqslant s_m(r)\leqslant \frac{1-y^{r+1}}{1-y}\binom{m}{r}$. These bounds are close if $\frac{r}{m}$ is near 0.
If $\frac{r}{m}$ is near $\frac{1}{2}$ then better approximations
involve the Berry-Esseen inequality~\cite{NC} to estimate
the binomial cumulative distribution function $2^{-m}s_m(r)$.
The cumulative distribution function
$\Phi(x)=\frac{1}{\sqrt{2\pi}}\int_{-\infty}^xe^{-t^2/2}\,dt$
is used in Remark~\ref{BE} to show that 
$|2^{-m}s_m(r)-\Phi(\frac{2r-m}{\sqrt{m}})|\leqslant\frac{0.4215}{\sqrt{m}}$ for
$0\leqslant r\leqslant m$ and $m\ne0$. Each binomial $\binom{m}{i}$ can
be estimated using Stirling's approximation as in~\cite[p.\,2]{Stanica}: 
$\binom{m}{i}=\frac{C_i^m}{\sqrt{2\pi p(1-p)m}}\left(1+O(\frac{1}{m})\right)$
where $C_i=\frac{1}{p^p(1-p)^{1-p}}$ and $p=p_i=i/m$. However, the sum
$\sum_{i=0}^r\binom{m}{i}$ of binomials is harder to approximate.
The preprint~\cite{Worsch} discusses different approximations to $s_m(r)$.

Sums of binomial coefficients modulo prime powers, where $i$ lies in a
congruence class, can be studied using number theory,
see~\cite[p.\,257]{Granville}.
Theorem~\ref{T:K} below shows how to find excellent rational approximations to
$s_m(r)$ via generalized continued fractions.

\begin{theorem}\label{T:K}
  Fix $r,m\in\Z$ where $0\leqslant r\leqslant m$ and recall that
  $s_m(r)=\sum_{i=0}^r\binom{m}{i}$.
  \begin{enumerate}[{\rm (a)}]
  \item If $b_i=2i(r+1-i)$, $a_i=m-2r+3i$ for $0\leqslant i\leqslant r$, then 
  \[
  Q\coloneq\frac{(r+1)\binom{m}{r+1}}{s_m(r)}=a_0+\Kop_{i=1}^{r} \frac{b_i}{a_i}.
  \]
  \item If $1\leqslant j\leqslant r$, then $\cT_j=R_j/R_{j-1}>0$ where
    $R_j\coloneq 2^j j!\sum_{k=0}^{r-j}\binom{r-k}{j}\binom{m}{k}$ satisfies
    $b_j R_{j-1} - a_j R_j = R_{j+1}$. Also, $(m-r)\binom{m}{r}-a_0R_0=R_1$.
  \end{enumerate}
\end{theorem}

Since $s_m(m)=2^m$, it follows that $s_m(m-r)=2^m-s_m(r-1)$ so we
focus on values of $r$ satisfying $0\leqslant r\leqslant\lfloor\frac{m}{2}\rfloor$.
Theorem~\ref{T:K}  allows us to find a sequence of successively sharper
upper and lower bounds for $Q$ (which can be made arbitrarily tight),
the coarsest being
$m-2r\leqslant Q\leqslant m-2r+\frac{2r}{m-2r+3}$ for $1\leqslant r<\frac{m+3}{2}$, see Proposition~\ref{P} and Corollary~\ref{C:approx}.

The fact that the tails $\cT_1,\dots,\cT_r$ are all positive is unexpected
as $b_i/a_i$ is negative if $\frac{m+3i}{2}<r$.
This fact is crucial for 
approximating $\cT_1=\Kop_{i=1}^{r} \frac{b_i}{a_i}$, see Theorem~\ref{T:root3}.
Theorem~\ref{T:K} implies that $\cT_1\cT_2\cdots\cT_r=R_r/R_0$.
Since $R_0=s_m(r)$, $R_r=2^rr!$, $\cT_j=\frac{b_j}{a_j+\cT_{j+1}}$ and
$\prod_{j=1}^rb_j=2^r(r!)^2$, the surprising~factorizations below
follow \emph{c.f.} Remark~\ref{R:fact}.

\begin{corollary}\label{C:beautiful}
  We have $s_m(r)\prod_{j=1}^r\cT_j=2^rr!$ and
  $r!s_m(r)=\prod_{j=1}^r(a_j+\cT_{j+1})$.
\end{corollary}

Suppose that $\omega>1$ and write $g(r)=g_{\omega,m}(r)$. We extend the domain
of $g(r)$ by setting $g(-1)=0$ and $g(m+1)=\frac{g(m)}{\omega}$ in
keeping with~\eqref{E1}. It is easy to prove that $g(r)$ is a
\emph{unimodal} function \emph{c.f.}~\cite[\S2]{BP}. Hence 
there exists some $r_0\in\{0,1,\dots,m\}$ that satisfies
\begin{equation}\label{E:unimodal}
  g_{\omega,m}(-1)<\cdots<g_{\omega,m}(r_0-1)\leqslant g_{\omega,m}(r_0)
  \quad\textup{and}\quad
  g_{\omega,m}(r_0)>\cdots>g_{\omega,m}(m+1).
\end{equation}
As $g(-1)<g(0)=1$ and $(\frac{2}{\omega})^m=g(m)>g(m+1)=\frac{2^m}{\omega^{m+1}}$, both chains of inequalities are non-empty.
The chains of inequalities~\eqref{E:unimodal} serve to define $r_0$.

We use Theorem~\ref{T:K} to show that $r_0$ is commonly close
to $r'\coloneq\lfloor\frac{m+2}{\omega+1}\rfloor$.
We always have $r'\leqslant r_0$ (by Lemma~\ref{L:floor}) and though $r_0-r'$
approaches $\frac{m}{2}$ as $\omega$ approaches 1
(see Remark~\ref{R:gap}), if $\omega\geqslant\sqrt{3}$ then
$0\leqslant r_0 - r'\leqslant 1$ by the next theorem.

\begin{theorem}\label{T:root3}
  If $\omega\geqslant\sqrt{3}$, $m\in\{0,1,\dots\}$
  and $r'\coloneq\lfloor\frac{m+2}{\omega+1}\rfloor$, then
  $r_0\in\{r',r'+1\}$, that is
  \[
    g(0)<\cdots<g(r'-1)\leqslant g(r'),\quad\textup{and}\quad
    g(r'+1)>g(r'+2)>\cdots>g(m).
  \]
\end{theorem}

Sharp bounds for $Q$ seem powerful: they enable
short and elementary proofs of results that previously required
substantial effort.  For example, 
our proof in~\cite[Theorem~1.1]{GP} for $\omega=2$ of the formula
$r_0=\lfloor\frac{m}{3}\rfloor+1$ involved a lengthy
argument, and our first proof of Theorem~\ref{T:formula}
below involved a delicate induction.
By this theorem there is a unique maximum, namely $r_0=r'=\lfloor\frac{m+2}{\omega+1}\rfloor$
when $\omega\in\{3,4,5,\dots\}$ and $\omega\ne m+1$, \emph{c.f.}~Remark~\ref{R:LLL}.
In particular, strict inequality $g_{\omega,m}(r'-1)<g_{\omega,m}(r')$ holds.

\begin{theorem}\label{T:formula}
  Suppose that $\omega\in\{3,4,5,\dots\}$ and
  $r'=\lfloor\frac{m+2}{\omega+1}\rfloor$. Then
  \[
  g_{\omega,m}(0)<\cdots<g_{\omega,m}(r'-1)\leqslant g_{\omega,m}(r')>g_{\omega,m}(r'+1)>\cdots>g_{\omega,m}(m),
  \]
  with equality if and only if $\omega=m+1$.
\end{theorem}

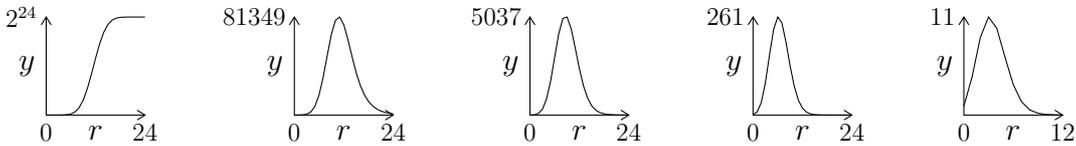
\begin{figure}[!ht]
\begin{center}
\begin{tikzpicture}[scale=1.3] 
\node [below,scale=0.8] (0,0) {$0$};
\node [left,scale=0.8] at (0,1) {$2^{24}$};
\node [below,scale=0.8] at (1,0) {$24$};
\node [left] at (0,0.5) {$y$};
\node [below] at (0.5,0) {$r$};
\draw [-] (0,1)-- (0,0) -- (1,0);
\draw [-] (0.92,0.05) -- (1,0) -- (0.92,-0.05);
\draw [-] (-0.05,0.92) -- (0,1) -- (0.05,0.92);
\draw 
(0.0000,0.0000)--
(0.04167,0.0000)--
(0.08334,0.0000)--
(0.1250,0.0001386)--
(0.1667,0.0007719)--
(0.2083,0.003305)--
(0.2500,0.01133)--
(0.2917,0.03196)--
(0.3333,0.07580)--
(0.3750,0.1537)--
(0.4167,0.2706)--
(0.4583,0.4194)--
(0.5000,0.5806)--
(0.5417,0.7294)--
(0.5833,0.8463)--
(0.6250,0.9242)--
(0.6667,0.9680)--
(0.7083,0.9886)--
(0.7500,0.9967)--
(0.7917,0.9992)--
(0.8333,0.9999)--
(0.8750,1.000)--
(0.9167,1.000)--
(0.9583,1.000)--
(1.000,1.000);
\end{tikzpicture}
\hskip5mm
\begin{tikzpicture}[scale=1.3] 
\node [below,scale=0.8] (0,0) {$0$};
\node [left,scale=0.8] at (0,1) {$81349$};
\node [below,scale=0.8] at (1,0) {$24$};
\node [left] at (0,0.5) {$y$};
\node [below] at (0.5,0) {$r$};
\draw [-] (0,1)-- (0,0) -- (1,0);
\draw [-] (0.92,0.05) -- (1,0) -- (0.92,-0.05);
\draw [-] (-0.05,0.92) -- (0,1) -- (0.05,0.92);
\draw 
(0.0000,0.0000)--
(0.04167,0.0002049)--
(0.08334,0.001644)--
(0.1250,0.008468)--
(0.1667,0.03145)--
(0.2083,0.08977)--
(0.2500,0.2051)--
(0.2917,0.3857)--
(0.3333,0.6099)--
(0.3750,0.8246)--
(0.4167,0.9678)--
(0.4583,1.000)--
(0.5000,0.9228)--
(0.5417,0.7729)--
(0.5833,0.5978)--
(0.6250,0.4352)--
(0.6667,0.3039)--
(0.7083,0.2070)--
(0.7500,0.1391)--
(0.7917,0.09296)--
(0.8333,0.06201)--
(0.8750,0.04134)--
(0.9167,0.02757)--
(0.9583,0.01838)--
(1.000,0.01225);
\end{tikzpicture}
\hskip5mm
\begin{tikzpicture}[scale=1.3] 
\node [below,scale=0.8] (0,0) {$0$};
\node [left,scale=0.8] at (0,1) {$5037$};
\node [below,scale=0.8] at (1,0) {$24$};
\node [left] at (0,0.5) {$y$};
\node [below] at (0.5,0) {$r$};
\draw [-] (0,1)-- (0,0) -- (1,0);
\draw [-] (0.92,0.05) -- (1,0) -- (0.92,-0.05);
\draw [-] (-0.05,0.92) -- (0,1) -- (0.05,0.92);
\draw 
(0.0000,0.0001985)--
(0.04167,0.002481)--
(0.08334,0.01494)--
(0.1250,0.05769)--
(0.1667,0.1607)--
(0.2083,0.3440)--
(0.2500,0.5895)--
(0.2917,0.8315)--
(0.3333,0.9861)--
(0.3750,1.000)--
(0.4167,0.8802)--
(0.4583,0.6820)--
(0.5000,0.4720)--
(0.5417,0.2965)--
(0.5833,0.1720)--
(0.6250,0.09393)--
(0.6667,0.04919)--
(0.7083,0.02512)--
(0.7500,0.01266)--
(0.7917,0.006347)--
(0.8333,0.003176)--
(0.8750,0.001588)--
(0.9167,0.0007941)--
(0.9583,0.0003970)--
(1.000,0.0001985);
\end{tikzpicture}
\hskip5mm
\begin{tikzpicture}[scale=1.3] 
\node [below,scale=0.8] (0,0) {$0$};
\node [left,scale=0.8] at (0,1) {$261$};
\node [below,scale=0.8] at (1,0) {$24$};
\node [left] at (0,0.5) {$y$};
\node [below] at (0.5,0) {$r$};
\draw [-] (0,1)-- (0,0) -- (1,0);
\draw [-] (0.92,0.05) -- (1,0) -- (0.92,-0.05);
\draw [-] (-0.05,0.92) -- (0,1) -- (0.05,0.92);
\draw 
(0.0000,0.003836)--
(0.04167,0.03197)--
(0.08334,0.1283)--
(0.1250,0.3303)--
(0.1667,0.6133)--
(0.2083,0.8754)--
(0.2500,1.000)--
(0.2917,0.9404)--
(0.3333,0.7435)--
(0.3750,0.5026)--
(0.4167,0.2950)--
(0.4583,0.1524)--
(0.5000,0.07030)--
(0.5417,0.02944)--
(0.5833,0.01139)--
(0.6250,0.004145)--
(0.6667,0.001447)--
(0.7083,0.0004928)--
(0.7500,0.0001656)--
(0.7917,0.0000)--
(0.8333,0.0000)--
(0.8750,0.0000)--
(0.9167,0.0000)--
(0.9583,0.0000)--
(1.000,0.0000);
\end{tikzpicture}
\hskip5mm
\begin{tikzpicture}[scale=1.3] 
\node [below,scale=0.8] (0,0) {$0$};
\node [below] at (0.5,0) {$r$};
\node [below,scale=0.8] at (1,0) {$12$};
\node [left,scale=0.8] at (0,1) {$11$};
\node [left] at (0,0.5) {$y$};
\draw [-] (0,1)-- (0,0) -- (1,0);
\draw [-] (0.92,0.05) -- (1,0) -- (0.92,-0.05);
\draw [-] (-0.05,0.92) -- (0,1) -- (0.05,0.92);
\draw
(0.0000,0.09030)--
(0.08334,0.3913)--
(0.1667,0.7926)--
(0.2500,1.000)--
(0.3333,0.8852)--
(0.4167,0.5894)--
(0.5000,0.3109)--
(0.5833,0.1363)--
(0.6667,0.05226)--
(0.7500,0.01843)--
(0.8333,0.006244)--
(0.9167,0.002087)--
(1.000,0.0006959);
\end{tikzpicture}
\end{center}
\vskip-8mm
\caption{Plots of $y=g_{\omega,24}(r)$ for $0\leqslant r\leqslant 24$ with $\omega\in\{1,\frac32,2,3\}$, and $y=g_{3,12}(r)$}\label{F1}
\end{figure}

Our motivation was to analyze $g_{\omega,m}(r)$ by using estimates for $Q$ given
by the generalized continued fraction in Theorem~\ref{T:K}. This gives
tighter estimates than the method involving partial sums used in~\cite{GP}.
The plots of $y=g_{\omega,m}(r)$ for $0\leqslant r\leqslant m$ are highly asymmetrical if
$\omega-1$ and $m$ are small. However, if $m$ is large the plots exhibit an `approximate symmetry'
about the vertical line~$r=r_0$ (see Figure~\ref{F1}). Our observation
that $r_0$ is
close to $r'$ for many choices of $\omega$ was the starting point of our
research.

Byun and Poznanovi\'c~\cite[Theorem~1.1]{BP} compute
the maximizing input, call it $r^*$, for the function
$f_{m,a}(r)\coloneq(1+a)^{-r}\sum_{i=0}^r\binom{m}{i}a^i$
where $a\in\{1,2,\dots\}$. Their function equals $g_{\omega,m}(r)$ only when
$\omega=1+a=2$. Some of their results and methods are similar to those
in~\cite{GP} which studied the case $\omega=2$. They prove that
$r^*=\lfloor\frac{a(m+1)+2}{2a+1}\rfloor$ provided
$m\not\in\{3,2a+4,4a+5\}$ or $(a,m)\ne(1,12)$ when
$r^*=\lfloor\frac{a(m+1)+2}{2a+1}\rfloor-1$.

In Section~\ref{S:gcf} we prove Theorem~\ref{T:K} and record approximations
to our generalized continued fraction expansion.
When $m$ is large, the plots of $y=g_{\omega,m}(r)$ are reminiscent of a normal
distribution with mean $\mu\approx\frac{m}{\omega+1}$. Section~\ref{S:max}
proves key lemmas for estimating~$r_0$,
and applies Theorem~\ref{T:K} to prove Theorem~\ref{T:formula}.
Non-integral values of $\omega$ are considered in Section~\ref{S:nonintegral}
where Theorem~\ref{T:root3} is proved.
In Section~\ref{S:est} we estimate the maximum height $g(r_0)$ using
elementary methods and estimations, see Lemma~\ref{L:bounds}.
A `statistical' approximation to $s_m(r)$ is given in Remark~\ref{BE}, and it
is compared in Remark~\ref{R:Approx} to the `generalized continued fraction
approximations' of $s_m(r)$ in Proposition~\ref{P}.

\section{Generalized continued fraction formulas}\label{S:gcf}

In this section we prove Theorem~\ref{T:K}, namely that
$Q\coloneq\frac{r+1}{s_m(r)}\binom{m}{r+1}=a_0+\cT_1$ where
$\cT_1=\Kop_{i=1}^r\frac{b_i}{a_i}$. The equality
$s_m(r)=\frac{r+1}{a_0+\cT_1}\binom{m}{r+1}$ is noted in
Corollary~\ref{C:GCF}.

A version of Theorem~\ref{T:K}(a) was announced in the SCS2022 Poster room, created to run concurrently with vICM 2022, see~\cite{P}.

\begin{proof}[Proof of Theorem~\ref{T:K}]
  
  Set $R_{-1}=Q\,s_m(r)=(r+1)\binom{m}{r+1}=(m-r)\binom{m}{r}$ and
  \[
  R_j=2^j j!\sum_{k=0}^{r-j}\binom{r-k}{j}\binom{m}{k}
  \qquad\textup{for $0\leqslant j\leqslant r+1$.}
  \]
  Clearly $R_0=s_m(r)$, $R_{r+1}=0$ and $R_j>0$ for $0\leqslant j\leqslant r$.
  We will prove in the following paragraph that the quantities $R_j$,
  $a_j=m-2r+3j$, and $b_j=2j(r+1-j)$
  satisfy the following $r+1$ equations, where the
  first equation~\eqref{E:first} is atypical:
  \begin{align}
    R_{-1}-a_0R_0&=R_1,\label{E:first}\\
    b_jR_{j-1}-a_jR_j&=R_{j+1}\quad\textup{where $1\leqslant j\leqslant r$.}\label{E:second}
  \end{align}
  Assuming~\eqref{E:second} is true, we prove by induction that
  $\cT_{j}=R_{j}/R_{j-1}$ holds for $r+1\geqslant j\geqslant 1$.
  This is clear for $j=r+1$ since $\cT_{r+1}=R_{r+1}=0$. Suppose that
  $1\leqslant j\leqslant r$ and $\cT_{j+1}=R_{j+1}/R_{j}$ holds. We show that
  $\cT_{j}=R_{j}/R_{j-1}$ holds. Using~\eqref{E:second} and $R_j>0$
  we have $b_jR_{j-1}/R_j-a_j=R_{j+1}/R_j=\cT_{j+1}$.
  Hence $R_j/R_{j-1}=b_j/(a_j+\cT_{j+1})=\cT_j$, completing the induction.
  Equation~\eqref{E:first} gives $Q=R_{-1}/R_0=a_0+R_1/R_{0}=a_0+\cT_1$ as claimed.
  Since $R_j>0$ for $0\leqslant j\leqslant r$, we have $\cT_{j}=R_{j}/R_{j-1}>0$ for $1\leqslant j\leqslant r$. This proves the first half of Theorem~\ref{T:K}(b), and the recurrence
  $\cT_j=b_j/(a_j+\cT_{j+1})$ for $1\leqslant j\leqslant r$, proves part~(a).

  We now show that~\eqref{E:first} holds.
  The identity $R_0=2^00!\sum_{k=0}^r\binom{m}{k}=s_m(r)$ gives
  \begin{align*}
    R_{-1}-a_0R_0&=(r+1)\binom{m}{r+1}-(m-2r)\sum_{i=0}^r\binom{m}{i}\\
    &=(r+1)\binom{m}{r+1}-\sum_{i=0}^r(-i+m-i-2r+2i)\binom{m}{i}\\
    &=\sum_{i=0}^{r} \left[(i+1)\binom{m}{i+1}-(m-i)\binom{m}{i}\right] +2\sum_{i=0}^{r-1}(r-i)\binom{m}{i}.
  \end{align*}
  As $(i+1)\binom{m}{i+1}=(m-i)\binom{m}{i}$, we get
  $R_{-1}-a_0R_0=2\sum_{k=0}^{r-1}\binom{r-k}{1}\binom{m}{k}=R_1$.
   
  We next show that~\eqref{E:second} holds.
  To simplify our calculations, we divide by $C_{j}\coloneq 2^j j!$.
  Using $(j+1)\binom{r-k}{j+1}=(r-k-j)\binom{r-k}{j}$ gives
  \begin{align*}
    \frac{R_{j+1}}{C_j}
    &=\sum_{k=0}^{r-j-1}2(j+1)\binom{r-k}{j+1}\binom{m}{k}\\
    &=\sum_{k=0}^{r-j}2(r-k-j)\binom{r-k}{j}\binom{m}{k}\\
    &=\sum_{k=0}^{r-j+1}(j-k)\binom{r-k}{j}\binom{m}{k}
      -\sum_{k=0}^{r-j}(k-2r+3j)\binom{r-k}{j}\binom{m}{k}
  \end{align*}
  noting that the term with $k=r-j+1$ in the first sum is zero as $\binom{j-1}{j}=0$.
  Writing $L=\sum_{k=0}^{r-j}(k-2r+3j)\binom{r-k}{j}\binom{m}{k}$ and using the identity
  $j\binom{r-k}{j}=(r+1-j-k)\binom{r-k}{j-1}$~gives
  \begin{align*}
    \frac{R_{j+1}}{C_j}=&\sum_{k=0}^{r-j+1}
    \left[(r+1-j-k)\binom{r-k}{j-1}-k\binom{r-k}{j}\right]\binom{m}{k}\; -L\\
    =&\sum_{k=0}^{r-j+1}
     \left[(r+1-j)\binom{r-k}{j-1}-k\binom{r-k}{j-1}-k\binom{r-k}{j}\right]\binom{m}{k} \;-L\\
    =&\sum_{k=0}^{r-j+1}
     \left[(r+1-j)\binom{r-k}{j-1}-k\binom{r-k+1}{j}\right]\binom{m}{k} \;-L.\\
  \end{align*}
  However, $k\binom{m}{k}=(m-k+1)\binom{m}{k-1}$, and therefore, 
  \begin{align*}
  \sum_{k=0}^{r-j+1}k\binom{r-k+1}{j}\binom{m}{k}&=\sum_{k=1}^{r-j+1}(m-k+1)\binom{r-k+1}{j}\binom{m}{k-1}\\&=\sum_{\ell=0}^{r-j}(m-\ell)\binom{r-\ell}{j}\binom{m}{\ell}.
  \end{align*}
  Thus
  \begin{align*}
    \frac{R_{j+1}}{C_{j}}
    &=\sum_{k=0}^{r-j+1}(r-j+1)\binom{r-k}{j-1}\binom{m}{k}
    -\sum_{k=0}^{r-j}(m-k)\binom{r-k}{j}\binom{m}{k} \;-L\\
    &=\sum_{k=0}^{r-j+1}(r-j+1)\binom{r-k}{j-1}\binom{m}{k}
    -\sum_{k=0}^{r-j}(m-k+k-2r+3j)\binom{r-k}{j}\binom{m}{k}\\
    &=\sum_{k=0}^{r-j+1}(r-j+1)\binom{r-k}{j-1}\binom{m}{k}
      - \sum_{k=0}^{r-j}\overbrace{(m-2r+3j)}^{a_j}\binom{r-k}{j}\binom{m}{k}\\
    &=\frac{\overbrace{2j(r-j+1)}^{b_j}2^{j-1}(j-1)!}{C_{j}}\sum_{k=0}^{r-j+1}\binom{r-k}{j-1}\binom{m}{k}
    - \sum_{k=0}^{r-j}a_j\binom{r-k}{j}\binom{m}{k}
  \end{align*}
  Hence $\frac{R_{j+1}}{C_{j}}=\frac{b_{j}R_{j-1}}{C_{j}}-\frac{a_{j}R_{j}}{C_{j}}$
  for $1\leqslant j\leqslant r$. When $j=r$, our convention gives $R_{r+1}=0$. This proves
  part~(b) and completes the proof of part~(a).
\end{proof}

\begin{remark}\label{R:fact}
  View $m$ as an indeterminant, so that $r!s_m(r)$
  is a polynomial in $m$ over $\Z$ of degree~$r$.
  The factorization $r!s_m(r)=\prod_{j=1}^r(a_j+\cT_{j+1})$ in
  Corollary~\ref{C:beautiful} involves the
  rational functions $a_j+\cT_{j+1}$. However, Theorem~\ref{T:K}(b) gives
  $\cT_{j+1}=\frac{R_{j+1}}{R_{j}}$, so that
  $a_j+\cT_{j+1}=\frac{a_{j}R_{j}+R_{j+1}}{R_j}=\frac{b_{j}R_{j-1}}{R_j}$.
  This determines the
  numerator and denominator of the rational function
  $a_j+\cT_{j+1}$, and explains why
  $\prod_{j=1}^r(a_j+\cT_{j+1})=
  \frac{R_0}{R_r}\prod_{j=1}^rb_{j}=r!s_m(r)$.
  This is different from, but reminiscent of, the ratio $p_{j+1}/p_j$
  described on p.\,26 of~\cite{Olds}\hfill$\diamond$
\end{remark}

\begin{corollary}\label{C:GCF}
  If $r,m\in\Z$ and  $0<r< m$, then
  \[
  s_m(r)\coloneq\sum_{i=0}^r\binom{m}{i}=\frac{(r+1)\binom{m}{r+1}}{m-2r+\cT_1}
  \qquad\textup{where $\cT_1=\Kop_{i=1}^r\frac{2i(r+1-i)}{m-2r+3i}>0$.}
  \]
\end{corollary}

If $r=0$, then $s_m(r)=\frac{(r+1)\binom{m}{r+1}}{m-2r+\cT_1}$ is true, but
$\cT_1=\Kop_{i=1}^r\frac{2i(r+1-i)}{m-2r+3i}=0$.

We will need some additional tools such as Proposition~\ref{P} and
Corollary~\ref{C:approx} below in order to prove Theorem~\ref{T:root3}.

Since $s_m(m-r)=2^m-s_m(r-1)$ approximating $s_m(r)$ for $0\leqslant r\leqslant m$
reduces to approximating $s_m(r)$ for
$0\leqslant r\leqslant \lfloor\frac{m}{2}\rfloor$. 
Hence the hypothesis $r<\frac{m+3}{2}$ in Proposition~\ref{P} and
Corollary~\ref{C:approx} is not too restrictive. Proposition~\ref{P}
generalizes~\cite[Theorem~3.3]{Olds}.

Let $\cH_j\coloneq\Kop_{i=1}^j\frac{b_i}{a_i}$ denote the \emph{$j$th head}
of the fraction $\Kop_{i=1}^r\frac{b_i}{a_i}$, where $\cH_0=0$.

\begin{proposition}\label{P}
  Let $b_i=2i(r+1-i)$ and $a_i=m-2r+3i$ for $0\leqslant i\leqslant r$.
  If $r<\frac{m+3}{2}$, then 
  $a_0+\cH_r=\frac{(r+1)\binom{m}{r+1}}{s_m(r)}$ can be approximated using
  the following chain of inequalities
  \[
  a_0+\cH_0<a_0+\cH_2<\cdots<a_0+\cH_{2\lfloor r/2\rfloor}
  <a_0+\cH_{2\lfloor (r-1)/2\rfloor+1}<\cdots<a_0+\cH_3<a_0+\cH_1.
  \]
\end{proposition}

\begin{proof}
  Note that $r$ equals either $2\lfloor r/2\rfloor$ or
  $2\lfloor (r-1)/2\rfloor+1$, depending on its parity.
  
  We showed in the proof of Theorem~\ref{T:K} that
  $\frac{(r+1)\binom{m}{r+1}}{s_m(r)}=a_0+\cH_r=a_0+\Kop_{i=1}^r\frac{b_i}{a_i}$.
  Since $r<\frac{m+3}{2}$, we have $a_i>0$ and $b_i>0$
  for $1\leqslant i\leqslant r$ and hence $\frac{b_i}{a_i}>0$. A straightforward
  induction (which we omit)
  depending on the parity of $r$ proves  that
  $\cH_0<\cH_2<\cdots<\cH_{2\lfloor r/2\rfloor}<\cH_{2\lfloor (r-1)/2\rfloor+1}<\cdots<\cH_3<\cH_1$. For example, if $r=3$, then
  \[
  \cH_0=0<\cfrac{b_1}{a_1 +\cfrac{b_2}{a_2}}
  <
  \cfrac{b_1}{a_1 +\cfrac{b_2}{a_2 +\cfrac{b_3}{a_3}}}
  <  \cfrac{b_1}{a_1}=\cH_1.
  \]
  proves $\cH_0<\cH_2<\cH_3<\cH_1$ as the tails are positive. Adding $a_0$ proves the claim.
\end{proof}

In asking whether $g_{\omega,m}(r)$ is a unimodal function, it
is natural to consider the ratio $g_{\omega,m}(r+1)/g_{\omega,m}(r)$
of successive terms. This suggests defining
\begin{equation}\label{E2}
  t(r)=t_{m}(r)\coloneq
  \frac{s_m(r+1)}{s_m(r)}=1+\frac{\binom{m}{r+1}}{s_m(r)}=1+\frac{Q}{r+1}.
\end{equation}
We will prove in Lemma~\ref{L:approx} that $t(r)$ is a strictly decreasing
function that determines when $g_{\omega,m}(r)$ is increasing or decreasing,
and $t_{m}(r_0-1)\geqslant \omega>t_{m}(r_0)$ determines $r_0$.

\begin{corollary}\label{C:approx}
  We have $m-2r\leqslant\frac{(r+1)\binom{m}{r+1}}{s_m(r)}$ for $r\geqslant 0$, and
  $\frac{(r+1)\binom{m}{r+1}}{s_m(r)}\leqslant m-2r+\frac{2r}{m-2r+3}$
  for $0\leqslant r<\frac{m+3}{2}$.
  Hence  $\frac{m+2}{r+1}\leqslant t_m(r)+1$ for $r\geqslant0$, and
  \[
  \frac{m+2}{r+1}\leqslant 
  t_m(r)+1\leqslant\frac{m+2}{r+1}+\frac{2r}{(r+1)(m-2r+3)}
  \quad\textup{for $0\leqslant r<\frac{m+3}{2}$.}
  \]
  Also $\frac{m+2}{r+1}< t_m(r)+1$ for $r>0$, and the above upper bound is
  strict for $1< r<\frac{m+3}{2}$.
\end{corollary}

\begin{proof}
  We proved
  $Q=\frac{(r+1)\binom{m}{r}}{s_m(r)}=(m-2r)+\Kop_{i=1}^r\frac{2i(r+1-i)}{m-2r+3i}$
  in Theorem~\ref{T:K}.
  Hence $m-2r=\frac{(r+1)\binom{m}{r+1}}{s_m(r)}$ if $r=0$ and
  $m-2r<\frac{(r+1)\binom{m}{r+1}}{s_m(r)}$ if $1\leqslant r<\frac{m+3}{2}$ by
  Proposition~\ref{P}. Clearly $m-2r<0\le\frac{(r+1)\binom{m}{r+1}}{s_m(r)}$
  if $\frac{m+3}{2}\le r\le m$. Similarly
  $\frac{(r+1)\binom{m}{r+1}}{s_m(r)}= m-2r+\frac{2r}{m-2r+3}$ if $r=0,1$, and
  again Proposition~\ref{P} shows that
  $\frac{(r+1)\binom{m}{r+1}}{s_m(r)}< m-2r+\frac{2r}{m-2r+3}$
  if $1< r<\frac{m+3}{2}$.
  The remaining inequalities (and equalities) follow similarly since
  $t_m(r)+1=2+\frac{\binom{m}{r+1}}{s_m(r)}$ and
  $2+\frac{m-2r}{r+1}=\frac{m+2}{r+1}$.
\end{proof}

\section{Estimating the maximizing input \texorpdfstring{$r_0$}{}}\label{S:max}

Fix $\omega>1$. In this section we consider
the function $g(r)=g_{\omega,m}(r)$ given by~\eqref{E1}.
As seen in Table~\ref{T:data}, it is easy to compute $g(r)$ if $r$ is
near $0$ or $m$. For $m$ large and $r$ near $0$, we have `sub-exponential'
growth $g(r)\approx\frac{m^r}{r!\omega^r}$. Similarly for $r$ 
near $m$, we have exponential decay $g(r)\approx\frac{2^m}{\omega^r}$.
The middle values require more thought.

\begin{table}[!ht]
  \vskip-3mm
  \caption{Values of $g_{w,m}(r)$}\label{T:data}
\begin{center}\vskip-3mm
  \begin{tabular}{rccccccccc} 
    \toprule
    $r$        &&0&              1&2&3&$\cdots$&$m-2$&$m-1$&$m$\\[2mm]
    $g_{\omega,m}(r)$&&1&$\frac{m+1}{\omega}$&$\frac{m^2+m+2}{2\omega^2}$&$\frac{m^3+5m+6}{6\omega^3}$&$\cdots$&$\frac{2^m-m-1}{\omega^{m-2}}$&$\frac{2^m-1}{\omega^{m-1}}$&$\left(\frac{2}{\omega}\right)^m$\\
    \bottomrule
\end{tabular}\label{T:tab}
\end{center}
\vskip-3mm
\end{table}

On the other hand, the plots $y=g(r)$, $0\leqslant r\leqslant m$, exhibit a remarkable
visual symmetry when $m$ is large. The relation $s_m(m-r)=2^m-s_m(r-1)$
and the distorting scale factor of $\omega^{-r}$ shape the plots. The 
examples in Figure~\ref{F1} show an approximate left-right symmetry about a
maximizing input $r\approx\frac{m}{\omega+1}$. It surprised the
authors that in many cases there exists a simple exact formula for the
maximizing input (it is usually unique as Corollary~\ref{C} suggests). In
Figure~\ref{F1} we have used different scale factors for the $y$-axes.
The maximum value of $g_{\omega,m}(r)$ varies considerably as $\omega$ varies
(\emph{c.f.} Lemma~\ref{L:bounds}), so we scaled
the maxima (rounded to the nearest integer) to the same height.

\begin{lemma}\label{L:approx}
  Recall that $g(r)=\omega^{-r} s_m(r)$ by~\eqref{E1} and
  $t(r)=\frac{s_m(r+1)}{s_m(r)}$ by~\eqref{E2}.
  \begin{enumerate}[{\rm (a)}]
  \item $t(r-1)>t(r)>\frac{m-r}{r+1}$ for $0\leqslant r\leqslant m$
    where $t(-1)\coloneq\infty$;
  \item $g(r)< g(r+1)$ if and only if $t(r)>\omega$;
  \item $g(r)\leqslant g(r+1)$ if and only if $t(r)\geqslant\omega$;
  \item $g(r)> g(r+1)$ if and only if $\omega> t(r)$;
  \item $g(r)\geqslant g(r+1)$ if and only if $\omega\geqslant t(r)$;
  \item if $\omega>1$ then some $r_0\in\{0,\dots,m\}$ satisfies
    $t(r_0-1)\geqslant \omega> t(r_0)$,  and this condition is equivalent to
    $g(0)<\cdots<g(r_0-1)\leqslant g(r_0)$ and $g(r_0)>\cdots>g(m)$.
  \end{enumerate}
\end{lemma}

\begin{proof}
  (a)~We prove, using induction on $r$, that
  $t(r-1)>t(r)>\binom{m}{r+1}/\binom{m}{r}$ holds for $0\leqslant r\leqslant m$.
  These inequalities are clear for $r=0$ as
  $\infty>m+1>m$. For real numbers $\alpha,\beta,\gamma,\delta>0$,
  we have $\alpha\delta-\beta\gamma>0$ if and only if
  $\frac{\alpha}{\beta}>\frac{\alpha+\gamma}{\beta+\delta}>\frac{\gamma}{\delta}$; that is, the mediant $\frac{\alpha+\gamma}{\beta+\delta}$
  of $\frac{\alpha}{\beta}$ and $\frac{\gamma}{\delta}$
  lies strictly between $\frac{\alpha}{\beta}$ and $\frac{\gamma}{\delta}$. If $0<r\leqslant m$, then by induction
  \[
    t(r-1)>\frac{\binom{m}{r}}{\binom{m}{r-1}}
    =\frac{m-r+1}{r}>\frac{m-r}{r+1}=\frac{\binom{m}{r+1}}{\binom{m}{r}}.
  \]
  Applying the `mediant sum' to $t(r-1)=\frac{s_m(r)}{s_m(r-1)}>\frac{\binom{m}{r+1}}{\binom{m}{r}}$ gives
  \[
  \frac{s_m(r)}{s_m(r-1)}>\frac{s_m(r)+\binom{m}{r+1}}{s_m(r-1)+\binom{m}{r}}
  =\frac{s_m(r+1)}{s_m(r)}=t(r)
  >\frac{\binom{m}{r+1}}{\binom{m}{r}}.
  \]
  Therefore $t(r-1)>t(r)>\binom{m}{r+1}/\binom{m}{r}=\frac{m-r}{r+1}$ completing the induction, and proving~(a).

  (b,c,d,e)~The following are equivalent: $g(r)< g(r+1)$;
  $\omega s_m(r)<s_m(r+1)$; and $\omega< t(r)$.
  The other claims are proved similarly by replacing $<$ with $\leqslant$, $>$, $\geqslant$.

  (f)~Observe that $t(m)=\frac{s_m(m+1)}{s_m(m)}=\frac{2^m}{2^m}=1$.
  By part~(a), the function $y=t(r)$  is decreasing for $-1\leqslant r\leqslant m$.  Since
  $\omega>1$, there exists an integer $r_0\in\{0,\dots,m\}$ such that
  $\infty=t(-1)>\cdots>t(r_0-1)\geqslant\omega>t(r_0)>\cdots>t(m)=1$.
  By parts (b,c,d,e) an equivalent condition is
  $g(0)<\cdots<g(r_0-1)\leqslant g(r_0)$ and $g(r_0)>\cdots>g(m)$.
\end{proof}

The following is an immediate corollary of Lemma~\ref{L:approx}(f).

\begin{corollary}\label{C}
  If $t(r_0-1)>\omega$, then the function $g(r)$ in~\eqref{E1} has a unique
  maximum at~$r_0$. If $t(r_0-1)=\omega$, then $g(r)$ has two equal maxima, one
  at $r_0-1$ and one at $r_0$.
\end{corollary}

As an application of Theorem~\ref{T:K} we show that the largest maximizing
input $r_0$ for $g_{\omega,m}(r)$ satisfies
$\lfloor\frac{m+2}{\omega+1}\rfloor\leqslant r_0$.
There are at most two maximizing inputs by Corollary~\ref{C}.

\begin{lemma}\label{L:floor}
  Suppose that $\omega>1$ and $m\in\Z$, $m\geqslant 0$.
  If $r'\coloneq\lfloor\frac{m+2}{\omega+1}\rfloor$, then
  $g(-1)<g(0)<\dots<g(r'-1)\leqslant g(r')$, and $g(-1)<\dots<g(r'-1)< g(r')$ if $r'>1$
  or $\omega\ne m+1$. 
\end{lemma}

\begin{proof}
  The result is clear when $r'=0$.
  If $r'=1$, then
  $r'\leqslant\frac{m+2}{\omega+1}$ gives $\omega\leqslant m+1$ or $g(0)\leqslant g(1)$.
  Hence $g(0)< g(1)$ if $\omega\ne m+1$.  Suppose that $r'>1$.
  By Lemma~\ref{L:approx}(c,f) the chain $g(0)<\dots< g(r')$ is
  equivalent to $g(r'-1)< g(r')$, that is $t(r'-1)> \omega$.
  However,
  $t(r'-1)+1>\frac{m+2}{r'}$ by Corollary~\ref{C:approx}
  and $r'\leqslant\frac{m+2}{\omega+1}$ implies
  $\frac{m+2}{r'}\geqslant\omega+1$. 
  Hence $t(r'-1)+1>\omega+1$, so that
  $t(r'-1)> \omega$ as desired.
\end{proof}

\begin{proof}[Proof of~Theorem~\ref{T:formula}]
  Suppose that $\omega\in\{3,4,\dots\}$.
  Then $g(0)<\dots<g(r'-1)\leqslant g(r')$ by Lemma~\ref{L:floor}
  with strictness when $\omega\ne m+1$.
  If $\omega=m+1$, then $r'=\lfloor\frac{m+2}{\omega+1}\rfloor=1$
  and $g(0)=g(1)$ as claimed. It remains to show that
  $g(r')>g(r'+1)>\cdots>g(m)$. However, we need
  only prove that $g(r')> g(r'+1)$ by Lemma~\ref{L:approx}(f),
  or equivalently $\omega> t(r')$ by Lemma~\ref{L:approx}(d).

  Clearly $\omega\geqslant3$ implies
  $r'\leqslant \frac{m+2}{\omega+1}\leqslant\frac{m+2}{4}$.
  As $0\leqslant r'<\frac{m+3}{2}$, Corollary~\ref{C:approx} gives
  \[
  \frac{m+2}{r'+1}+\frac{2r'}{(r'+1)(m-2r'+3)}\geqslant t(r')+1.
  \]
  Hence $\omega+1> t(r')+1$ holds if
  $\omega+1> \frac{m+2}{r'+1}+\frac{2r'}{(r'+1)(m-2r'+3)}$.
  Since $\omega+1$ is an integer, we have $m+2=r'(\omega+1)+c$
  where $0\leqslant c\leqslant \omega$. It follows from
  $0\leqslant r'\leqslant\frac{m+2}{4}$
  that $\frac{2r'}{m-2r'+3}<1$.  This inequality
  and $m+2\leqslant r'(\omega+1)+\omega$ gives
  \[
    m+2+\frac{2r'}{m-2r'+3}<r'(\omega+1)+\omega+1=(r'+1)(\omega+1).
  \]
  Thus $\omega+1>\frac{m+2}{r'+1}+\frac{2r'}{(r'+1)(m-2r'+3)}\geqslant t(r')+1$,
  so $\omega>t(r')$ as required.
\end{proof}
  
\begin{remark}\label{R:omega2}
    The proof of Theorem~\ref{T:formula} can be adapted to the case $\omega=2$.
    If $m+2=3r'+c$ where $c\leqslant\omega-1=1$, then
    $\frac{2r'}{m-2r'+3}=\frac{2r'}{r'+c+1}<2$, and if
    $c=\omega=2$, then a sharper $\cH_2$-bound must be used.
    This leads to a much shorter proof than~\cite[Theorem~1.1]{GP}.\hfill$\diamond$
\end{remark}

\section{Non-integral values of \texorpdfstring{$\omega$}{}}\label{S:nonintegral}

In this section, we prove that the maximum value of $g(r)$ is $g(r')$
or $g(r'+1)$ if $\omega\geqslant\sqrt{3}$. Before proving this
result (Theorem~\ref{T:root3}), we shall prove two preliminary lemmas.

\begin{lemma}\label{L:Gerhard}
  Suppose that $\omega>1$ and $r'\coloneq\lfloor\frac{m+2}{\omega+1}\rfloor$.
  If $\frac{m+2}{r'+1}\geqslant\sqrt{3}+1$, then
  \[
    g(-1)<g(0)<\cdots<g(r'-1)\leqslant g(r'),\quad\textup{and}\quad
    g(r'+1)>g(r'+2)>\cdots>g(m).
  \]
\end{lemma}

\begin{proof}
  It suffices,
  by Lemma~\ref{L:approx}(f) and Lemma~\ref{L:floor}
  to prove that $g(r'+1)>g(r'+2)$.   The strategy is to
  show $\omega>t(r'+1)$, that is $\omega+1>t(r'+1)+1$. However,
  $\omega+1>\frac{m+2}{r'+1}$, so it suffices to prove
  that $\frac{m+2}{r'+1}\geqslant t(r'+1)+1$.
  Since $r'+1\leqslant\frac{m+2}{\sqrt{3}+1}<\frac{m+2}{2}$, we can use
  Corollary~\ref{C:approx} and just prove that
  $\frac{m+2}{r'+1}\geqslant\frac{m+2}{r'+2}+\frac{2r'+2}{(r'+2)(m-2r'+1)}$.
  This inequality is equivalent to
  $\frac{m+2}{r'+1}\geqslant\frac{2r'+2}{m-2r'+1}$.
  However, $\frac{m+2}{r'+1}\geqslant\sqrt{3}+1$, so we need only show that
  $\sqrt{3}+1\geqslant\frac{2(r'+1)}{m-2r'+1}$, or equivalently 
  $m-2r'+1\geqslant(\sqrt{3}-1)(r'+1)$. This is true since
  $\frac{m+2}{r'+1}\geqslant\sqrt{3}+1$ implies   
  $m-2r'+1\geqslant(\sqrt{3}-1)r'+\sqrt{3}>(\sqrt{3}-1)(r'+1)$.
\end{proof}

\begin{remark}\label{R:LLL}
  The strict inequality $g(r'-1)<g(r')$ holds by Lemma~\ref{L:floor}
  if $r'>1$ or $\omega\ne m+1$. It holds vacuously for $r'=0$.
  Hence adding the additional hypothesis that $\omega\ne m+1$ if $r'=1$
  to Lemma~\ref{L:Gerhard} (and Theorem~\ref{T:root3}),
  we may conclude that the inequality $g(r'-1)\leqslant g(r')$ is strict.\hfill$\diamond$
\end{remark}

\begin{remark}
  In Lemma~\ref{L:Gerhard}, the maximum can occur at $r'+1$.
  If $\omega=2.5$ and $m=8$, then $r'=\lfloor\frac{10}{3.5}\rfloor=2$ and
  $\frac{m+2}{r'+1}=\frac{10}{3}\geqslant\sqrt{3}+1$ however
  $g_{2.5,8}(2)=\frac{740}{125}<\frac{744}{125}=g_{2.5,8}(3)$.\hfill$\diamond$
\end{remark}

\begin{remark}\label{R:gap}
  The gap between $r'$ and the largest maximizing input $r_0$ can be
  arbitrarily large if $\omega$ is close to 1. For $\omega>1$, we have
  $r'=\lfloor\frac{m+2}{\omega+1}\rfloor<\frac{m+2}{2}$.
  If $1<\omega\leqslant\frac{1}{1-2^{-m}}$, then $g(m-1)\leqslant g(m)$, so $r_0=m$.
  Hence $r_0-r'>\frac{m-2}{2}$.\hfill$\diamond$
\end{remark}

\begin{remark}
  Since $r'\leqslant\lfloor\frac{m+2}{\omega+1}\rfloor<r'+1$, we see that
  $r'+1\approx\frac{m+2}{\omega+1}$, so that
  $\frac{m+2}{r'+1}\approx\omega+1$.
  Thus Lemma~\ref{L:Gerhard} suggests that if $\omega\gtrsim\sqrt{3}$, then $g_{\omega,m}(r)$ may have a
  maximum at $r'$ or $r'+1$. This heuristic reasoning is made
  rigorous in Theorem~\ref{T:root3}.\hfill$\diamond$
\end{remark}

\begin{remark}\label{R:K}
  Theorem~\ref{T:K} can be rephrased as
  $t_m(r)=\frac{s_m(r+1)}{s_m(r)}=\frac{m-r+1}{r+1}+\frac{\cK_m(r)}{r+1}$
  where
  \begin{equation}\label{E:Kr}
  \cK_m(r)\coloneq\Kop_{i=1}^{r} \frac{2i(r+1-i)}{m-2r+3i}=
  \cfrac{2r}{m-2r+3 +\cfrac{4r-4}{m-2r+6 +\cfrac{6r-12}{\llap{\ensuremath{\ddots}} \raisebox{-0.8em}{\ensuremath{m+r-3 + \cfrac{2r}{m+r}}}}}}.
  \end{equation}
  The following lemma repeatedly uses the expression $\omega>t_m(r+1)$.
  This is equivalent to
  $\omega>\frac{m-r}{r+2}+\frac{\cK_m(r+1)}{r+2}$, that is
  $(\omega+1)(r+2)>m+2+\cK_m(r+1)$.\hfill$\diamond$
\end{remark}

\begin{lemma}\label{L:root3}
  Let $m\in\{0,1,\dots\}$ and $r'=\lfloor\frac{m+2}{\omega+1}\rfloor$.
  If any of the following three conditions are met,
  then $g_{\omega,m}(r'+1)>\cdots>g_{\omega,m}(m)$ holds:
  \begin{center}
    {\rm (a)} $\omega\geqslant2$, or \quad
    {\rm(b)} $\omega\geqslant\frac{1+\sqrt{97}}{6}$ and $r'\ne2$, or \quad
    {\rm(c)} $\omega\geqslant\sqrt{3}$ and $r'\not\in\{2,3\}$.
  \end{center}
\end{lemma}

\begin{proof}
  The conclusion $g_{\omega,m}(r'+1)>\cdots>g_{\omega,m}(m)$ holds trivially if
  $r'+1\geqslant m$. Suppose henceforth that $r'+1< m$. Except for the excluded
  values of $r',\omega$, we will prove that
  $g_{\omega,m}(r'+1)>g_{\omega,m}(r'+2)$ holds, as this implies
  $g_{\omega,m}(r'+1)>\cdots>g_{\omega,m}(m)$ by Lemma~\ref{L:approx}(f).
  Hence we must prove that $\omega>t_m(r'+1)$ by Lemma~\ref{L:approx}(d).
  
  Recall that $r'\leqslant\frac{m+2}{\omega+1}<r'+1$. If $r'=0$, then
  $m+2<\omega+1$, that is $\omega>m+1>t(1)$ as desired. Suppose now that $r'=1$.
  There is nothing to prove if $m=r'+1=2$.
    Assume that $m>2$. Since $m+2<2(\omega+1)$, we have
    $2<m<2\omega$. The last
  line of Remark~\ref{R:K} and~\eqref{E:Kr} give the desired inequality:
  \[
    \omega>\frac{m}{2}\geqslant\frac{m-1}{3}+\frac{4}{3\left(m-1+\cfrac{4}{m+2}\right)}=t_m(2).
  \]
  In summary, $g_{\omega,m}(r'+1)>\cdots>g_{\omega,m}(m)$ holds for all $\omega>1$
  if $r'\in\{0,1\}$.

  We next prove $g_{\omega,m}(r'+1)>g_{\omega,m}(r'+2)$, or equivalently
  $\omega>t_m(r'+1)$ for $r'$ large enough, depending on $\omega$.
  We must prove that $(\omega+1)(r'+2)>m+2+\Kop_m(r'+1)$
  by Remark~\ref{R:K}.
  Writing $m+2=(\omega+1)(r'+\eps)$ where $0\leqslant\eps<1$, our goal, therefore,
  is to show $(\omega+1)(2-\eps)>\cK_m(r'+1)$.
  Using~\eqref{E:Kr} gives
  \[
  \cK_m(r'+1)=\frac{2(r'+1)}{m-2(r'+1)+3+\cT}
  =\frac{2(r'+1)}{(\omega+1)(r'+\eps)-2(r'+1)+1+\cT}
  \]
  where $\cT>0$ by Theorem~\ref{T:K} as $r'>0$.
  Rewriting the denominator using
  \[
  (\omega+1)(r'+\eps)-2(r'+1)=(\omega-1)(r'+1)-(\omega+1)(1-\eps),
  \]
  our goal $(\omega+1)(2-\eps)>\cK_m(r'+1)$ becomes
  \[
  (\omega+1)(2-\eps)\left[(\omega-1)(r'+1)-(\omega+1)(1-\eps)+1+\cT\right]
  >2(r'+1).
  \]
  Dividing by $(2-\eps)(r'+1)$ and rearranging gives
  \[
  (\omega^2-1)+\frac{(\omega+1)(1+\cT)}{r'+1}>\frac{2}{2-\eps}
  +\frac{(\omega+1)^2(1-\eps)}{r'+1}.
  \]
  This inequality may be written
  $(\omega^2-1)+\lambda>\frac{2}{2-\eps}+\mu(1-\eps)$ where
  $\lambda=\frac{(\omega+1)(1+\cT)}{r'+1}>0$ and
  $\mu=\frac{(\omega+1)^2}{r'+1}>0$. We view
  $f(\eps)\coloneq\frac{2}{2-\eps}+\mu(1-\eps)$ as a function of a real
  variable~$\eps$ where $0\leqslant\eps<1$. However, $f(\eps)$ is concave as the
  second derivative $f''(\eps)=\frac{4}{(2-\eps)^3}$ is positive for
  $0\leqslant\eps<1$.
  Hence the maximum value occurs at an end point: either $f(0)=1+\mu$ or
  $f(1)=2$. Therefore, it suffices to prove that
  $(\omega^2-1)+\lambda>\max\{2,1+\mu\}$.

  If $2\geqslant1+\mu$, then the desired bound $(\omega^2-3)+\lambda>0$ holds
  as $\omega\geqslant\sqrt{3}$. Suppose now that $2<1+\mu$. We must show
  $(\omega^2-1)+\lambda>1+\mu$, that is
  $\omega^2-2>\mu-\lambda=\frac{(\omega+1)(\omega-\cT)}{r'+1}$. Since
  $\cT>0$, a stronger inequality (that implies this)
  is $\omega^2-2\geqslant\frac{(\omega+1)\omega}{r'+1}$. The (equivalent) quadratic
  inequality $r'\omega^2-\omega-2(r'+1)\geqslant 0$ in $\omega$ is true provided
  $\omega\geqslant \frac{1+\sqrt{1+8r'(r'+1)}}{2r'}$. This says $\omega\geqslant2$
  if $r'=2$, and $\omega\geqslant\frac{1+\sqrt{97}}{6}$ if $r'=3$.
  If $r'\geqslant4$, we have
  \[
  \frac{1+\sqrt{1+8r'(r'+1)}}{2r'}=
  \frac{1}{2r'}+\sqrt{\frac{1}{4(r')^2}+2\left(1+\frac{1}{r'}\right)}
  \leqslant\frac{1}{8}+\sqrt{\frac{1}{64}+\frac{5}{2}}<\sqrt{3}.
  \]
  The conclusion now follows from the fact that
  $2>\frac{1+\sqrt{97}}{6} >\sqrt{3}$.
\end{proof}

\begin{proof}[Proof of Theorem~\ref{T:root3}]
  By Lemma~\ref{L:Gerhard} it suffices to show that $g(r'+1)>g(r'+2)$
  holds when $r'+1<m$ and $\omega\geqslant\sqrt{3}$.
  By Lemma~\ref{L:root3}(a), we can assume that $\sqrt{3}\leqslant\omega<2$ and
  $r'\in\{2,3\}$. For these choices of $\omega$ and $r'$, we must show
  that $\omega>t_m(r'+1)$ by Lemma~\ref{L:approx} for all permissible
  choices of $m$. Since $(\omega+1)r'\leqslant m+2<(\omega+1)(r'+1)$, when
  $r'=2$ we have $5<2(\sqrt{3}+1)\leqslant m+2<9$ so that $4\leqslant m\leqslant 6$.
  However, $t_m(3)$ equals $\frac{16}{15},\frac{31}{26},\frac{19}{14}$
  for these values of $m$. Thus $\sqrt{3}>t_m(3)$ holds as desired.
  Similarly, if $r'=3$, then $8<3(\sqrt{3}+1)\leqslant m+2<12$ so that $7\leqslant m\leqslant 9$.
  In this case $t_m(4)$ equals $\frac{40}{33},\frac{219}{163},\frac{191}{128}$
  for these values of $m$. In each case $\sqrt{3}>t_m(4)$, so the proof
  is complete.
\end{proof}

\begin{remark}\label{R:d}
  We place Remark~\ref{R:gap} in context.
  The conclusion of Theorem~\ref{T:root3} remains true for values of
  $\omega$ smaller than $\sqrt{3}$ and not `too close to~1'
  and $m$ is `sufficiently large'.
  Indeed, by adapting the proof of Lemma~\ref{L:root3} we can show
  there exists a sufficiently large integer $d$ such
  that $m > d^4$ and $\omega > 1 + \frac{1}{d}$ implies
  $g(r'+d) > g(r'+d+1)$. This shows that
  $r' \leqslant r_0 \leqslant r' +d$, so $r_0-r'\leqslant d$.
  We omit the technical proof of this fact.~$\hfill\diamond$
\end{remark}

\begin{remark}
  The sequence, $a_0+\cH_1,\dots,a_0+\cH_r$ terminates at
  $\frac{r+1}{s_m(r)}\binom{m}{r+1}$ by Theorem~\ref{T:K}.
  We will not comment here on \emph{how quickly} the alternating sequence
  in Proposition~\ref{P} converges when $r<\frac{m+3}{2}$.
  If $r=m$, then $a_0=-m$ and $\frac{m+1}{s_m(m+1)}\binom{m}{m+1}=0$,
  so Theorem~\ref{T:K} gives the curious  identity
  $\cH_m=\Kop_{i=1}^m\frac{2i(m+1-i)}{3i-m}=m$. If $\omega$ is less
  than $\sqrt{3}$ and `not too close to 1', then we believe that
  $r_0$ is approximately
  $\lfloor\frac{m+2}{\omega+1}+ \frac{2}{\omega^2 -1}\rfloor$,
  \emph{c.f.} Remark~\ref{R:d}.
\end{remark}

\section{Estimating the maximum value of \texorpdfstring{$g_{\omega,m}(r)$}{}}\label{S:est}

In this section we relate the size of the maximum value $g_{\omega,m}(r_0)$ to
the size of the binomial coefficient $\binom{m}{r_0}$.
In the case that we
know a formula for a maximizing input $r_0$, we can readily
estimate~$g_{\omega,m}(r_0)$ using approximations, such as~\cite{Stanica},
for binomial coefficients.

\begin{lemma}\label{L:bounds}
  The maximum value $g_{\omega,m}(r_0)$ of $g_{\omega,m}(r)$, $0\leqslant r\leqslant m$,
  satisfies
  \[
  \frac{1}{(\omega-1)\omega^{r_0}}\binom{m}{r_0+1}
  < g_{\omega,m}(r_0)
  \leqslant \frac{1}{(\omega-1)\omega^{r_0-1}}\binom{m}{r_0}.
  \]
\end{lemma}

\begin{proof}
  Since $g(r_0)$ is a maximum value, we
  have $g(r_0-1)\leqslant g(r_0)$. This is equivalent to
  $(\omega-1)s_m(r_0-1)\leqslant \binom{m}{r_0}$ as
  $s_m(r_0)=s_m(r_0-1)+\binom{m}{r_0}$.
  Adding $(\omega-1)\binom{m}{r_0}$ to both sides gives the equivalent
  inequality  $(\omega-1) s_m(r_0)\leqslant \omega\binom{m}{r_0}$. This proves
  the upper bound.
  
  Similar reasoning shows that the following are equivalent:
  (a)~$g(r_0)> g(r_0+1)$; (b)~$(\omega-1)s_m(r_0)> \binom{m}{r_0+1}$;
  and (c)~$g_{\omega,m}(r_0)>\frac{1}{(\omega-1)\omega^{r_0}}\binom{m}{r_0+1}$.
\end{proof}

In Theorem~\ref{T:root3} the maximizing input $r_0$ satisfies $r_0=r'+d$
where $d\in\{0,1\}$. In such cases when $r_0$ and $d$ are known,
we can bound the maximum $g_{\omega,m}(r_0)$ as follows.
  
\begin{corollary}
  Set $r' \coloneq \lfloor\frac{m+2}{\omega+1}\rfloor$ and
  $k \coloneq m+2 - (\omega+1)r'$.  Suppose that $r_0 = r'+d$ and
  $G=\frac{1}{(\omega-1)\omega^{r_0-1}}\binom{m}{r_0}$.
  Then $0\leqslant k< \omega+1$, $d\geqslant0$ and
  $1-\frac{1+d-\frac{d+2-k}{\omega}}{r_0+1}
    <\frac{g_{\omega,m}(r_0)}{G} \leqslant 1$.
\end{corollary}

\begin{proof}
  By Lemma~\ref{L:floor}, $r_0 = r'+d$ where $d\in\{0,1,\dots\}$. Since
  $r' = \lfloor\frac{m+2}{\omega+1}\rfloor$, we have
  $m+2 = (\omega+1)r'+k$ where $0\leqslant k< \omega+1$.
  The result follows from Lemma~\ref{L:bounds} and
  $m = (\omega+1)(r_0-d)+k-2$ as
  $\binom{m}{r_0+1}=\frac{m-r_0}{r_0+1}\binom{m}{r_0}$ and $\frac{m-r_0}{r_0+1}$ equals
  \[
  \frac{\omega(r_0-d)-d+k-2}{r_0+1}
  =\omega-\frac{\omega+\omega d+d+2-k}{r_0+1}
  =\omega\left(1-\frac{1+d+\frac{d+2-k}{\omega}}{r_0+1}\right).\qedhere
  \]
\end{proof}

The following remark is an application of the Chernoff bound,
\emph{c.f.}~\cite[Section~4]{Worsch}. Unlike Theorem~\ref{T:K}, it requires
the cumulative distribution function $\Phi(x)$, which is a
non-elementary integral, to approximate $s_m(r)$. It seems to give better
approximations only for values of $r$ near $\frac{m}{2}$,
see Remark~\ref{R:Approx}.

\begin{remark}\label{BE}
  We show how the Berry-Esseen inequality for a sum of binomial random
  variables can be used to approximate $s_m(r)$.
  Let $B_1,\dots,B_m$ be independent identically distributed binomial
  variables with a parameter $p$ where $0<p<1$, so that $P(B_i=1)=p$ and
  $P(B_i=0)=q\coloneq 1-p$. Let $X_i\coloneq B_i-p$ and
  $X\coloneq\frac{1}{\sqrt{mpq}}(\sum_{i=1}^m X_i)$. Then
  \[
  E(X_i)=E(B_i)-p=0,\quad E(X_i^2)=pq,\quad {\rm and}\quad
  E(|X_i|^3)=pq(p^2+q^2).
  \]
  Hence $E(X)=\frac{1}{\sqrt{mpq}}(\sum_{i=1}^m E(X_i))=0$
  and $E(X^2)=\frac{1}{mpq}(\sum_{i=1}^m E(X_i^2))=1$.
  By~\cite[Theorem~2]{NC} the Berry-Esseen inequality
  applied to $X$ states that
  \[
  |P(X\leqslant x)-\Phi(x)|\leqslant\frac{Cpq(p^2+q^2)}{(pq)^{3/2}\sqrt{m}}
  =\frac{C(p^2+q^2)}{\sqrt{mpq}}
  \qquad\textup{for all $m\in\{1,2,\dots\}$ and $x\in\R$,}
  \]
  where the constant $C\coloneq 0.4215$ is close to the
  lower bound $C_0=\frac{10+\sqrt{3}}{6\sqrt{2\pi}}=0.4097\cdots$
  and $\Phi(x)=\frac{1}{\sqrt{2\pi}}\int_{-\infty}^x e^{-t^2/2}\,dt=\frac{1}{2}(1+{\rm erf}(\frac{x}{\sqrt{2}}))$ is the cumulative distribution function for
  standard normal distribution.
  
  Writing $B=\sum_{i=1}^mB_i$ we have
  $P(B\leqslant b)=\sum_{i=0}^{\lfloor b\rfloor}\binom{m}{i}p^iq^{m-i}$ for $b\in\R$. 
  Thus $X=\frac{B-mp}{\sqrt{mpq}}$ and $x=\frac{b-mp}{\sqrt{mpq}}$ satisfy
  \[
    \left|P(B\leqslant b)-\Phi\left(\frac{b-mp}{\sqrt{mpq}}\right)\right|
    \leqslant\frac{C(p^2+q^2)}{\sqrt{mpq}}
    \qquad\textup{for all $m\in\{1,2,\dots\}$ and $b\in\R$.}
  \]
  Setting $p=q=\frac{1}{2}$, and taking $b=r\in\{0,1,\dots,m\}$ shows
  \[
    \left|2^{-m}s_m(r)-\Phi\left(\frac{2r-m}{\sqrt{m}}\right)\right|
    \leqslant\frac{0.4215}{\sqrt{m}}
    \qquad\textup{for $m\in\{1,2,\dots\}$.}
    \qquad\qquad\qquad\qquad\qquad\diamond
  \]
\end{remark}

\begin{remark}\label{R:Approx}
  Let $a_0+\cH_k$ be the generalized continued fraction approximation
  to $\frac{(r+1)\binom{m}{r+1}}{s_m(r)}$ suggested by Theorem~\ref{T:K}, where
  $\cH_k\coloneq\Kop_{i=1}^k\frac{b_i}{a_i}$, and $k$ is the depth of the
  generalized continued fraction. We compare the following two quantities:
  \[
  e_{m,r,k}\coloneq1-\frac{(r+1)\binom{m}{r+1}}{(a_0+\cH_k)s_m(r)}\qquad\textup{and}\qquad
  E_{m,r}\coloneq\left|1-\frac{2^m\Phi(\frac{2r-m}{\sqrt{m}})}{s_m(r)}\right|
  \leqslant\frac{0.4215\cdot 2^m}{\sqrt{m}\,s_m(r)}.
  \]
  The sign of $e_{m,r,k}$ is governed by the parity of $k$ by
  Proposition~\ref{P}.
  We shall assume that $r\leqslant\frac{m}{2}$.
  As $\frac{2^m}{s_m(r)}$ ranges from $2^m$ to about 2 as $r$ ranges from
  0 to $\lfloor\frac{m}{2}\rfloor$, it
  is clear that the upper bound for $E_{m,r}$ will be huge
  unless $r$ satisfies $\frac{m-\varepsilon}{2}\le r\le \frac{m}{2}$
  where $\varepsilon$ is `small' compared to $m$.
  By contrast, the computer code~\cite{G} verifies
  that the same is true for $E_{m,r}$, and shows that $|e_{m,r,k}|$ is
  small, even when $k$ is tiny, when $0\le r<\frac{m-\varepsilon}{2}$,
  see Table~\ref{T:E}.
  Hence the `generalized continued fraction' approximation to $s_m(r)$ is
  complementary to the `statistical' approximation, as shown in Table~\ref{T:E}.
  The reader can extend Table~\ref{T:E} by running the
  code~\cite{G} written in the {\sc Magma}~\cite{BosmaEtAl} language,
  using the online calculator \url{http://magma.maths.usyd.edu.au/calc/},
  for example.
  \vskip-4mm
  \begin{table}[!ht]
    \caption{Upper bounds for $|e_{m,r,k}|$ and $E_{m,r}$ for $m=10^4$}\label{T:E}
    \begin{center}
    \vskip-3mm
    \begin{tabular}{cllll}\toprule
      $r$&$|e_{m,r,3}|$&$|e_{m,r,5}|$&$|e_{m,r,21}|$&$E_{m,r}$\\
      1000&$2.3\times10^{-17}$&$6.6\times10^{-25}$&$5.7\times10^{-79}$&$1$\\
      4500&$1.3\times10^{-7}$&$2.5\times10^{-10}$&$7.1\times10^{-27}$&$0.018$\\
      5000&$0.93$&$0.86$&$0.24$&$0.008$\\ \bottomrule
      \end{tabular}
    \end{center}
  \end{table}
  \vskip-13mm\hfill$\diamond$
\end{remark}

\vskip-10mm
\paragraph*{Acknowledgment} 
SPG received support from the Australian Research Council Discovery
Project Grant DP190100450. GRP thanks his family, and SPG thanks his mother.

\end{document}